\documentclass[11pt,psamsfont]{amsart} \usepackage{bkmacros}
\include{graphicx}

\begin{document} 
%\begingroup\def\uppercasenonmath#1{}%
\let\MakeUppercase\relax%

\title{CAT(0) spaces with polynomial divergence of geodesics}
\author[Nata\v{s}a Macura]{Nata\v{s}a Macura\\ Department of
  Mathematics\\ Trinity University\\ nmacura@trinity.edu\\ One Trinity
  Place\\ San Antonio TX 78212}


\begin{abstract} We construct  a family of finite 2-complexes whose universal
  covers are CAT(0) and have polynomial divergence of desired
  degree. This answers a question of Gersten, namely whether such
  CAT(0) complexes exist.\end{abstract}

 \maketitle
\section {Introduction} 

In \cite{gerstendiv2} Gersten defined divergence of a CAT(0) space,
generalizing the classical idea of the divergence of geodesics in
manifolds, and showed it to be a quasi-isometry invariant. He
constructed a CAT(0) 2-complex with quadratic divergence, therefore
showing that the aphorism of Riemannian geometry that geodesics
diverge either linearly or exponentially fails for CAT(0) spaces. In
later work, Gersten \cite{gerstendiv2}, and M. Kapovich and
Leeb,\cite{mishadiv} showed that the aphorism also fails for
3-manifolds since there exist graph manifolds with quadratic
divergence of geodesics.  In this paper we exhibit a family of CAT(0)
groups $G_d, d\in \mathbb{N}$, such that the divergence of $G_d$ is
polynomial of degree $d.$ We construct $G_d$ inductively as an HNN
extension of $G_{d-1},$ starting with
$G_1=\mathbb{Z}\oplus\mathbb{Z}$. Each $G_d$ has a 2-dimensional
presentation complex $X_d$ whose universal cover $\wtd{X}_d$ is a
CAT(0) cube complex.  We prove that the divergence of $\wtd{X}_d$ is
polynomial of degree $d.$ The groups described here turn out to be the
family of examples W. Dison and T. Riley introduced and named hydra
groups in \cite{hydra}. W. Dison and T. Riley show that hydra groups
have finite-rank free subgroups with huge distortion and use this
class of groups to construct elementary examples of groups whose Dehn
functions are equally large.

Divergence of geodesics, as well as in its higher dimensional
generalizations received renewed interest in recent work of a number
of authors. In \cite{RAAGdivN} A. Abrams, N.  Brady, P. Dani,
M. Duchin and R. Young define higher divergence functions, which
measure isoperimetric properties "at infinity", and give a
characterization of the divergence of geodesics in RAAGs as well as
upper bound for filling loops at infinity in the mapping class group.
J. Behrstock and R. Charney (\cite{BehrChar}) give a group theoretic
characterization of geodesics with super-linear divergence in the
Cayley graph of a right-angled Artin group $A_{\Gamma}$ with connected
defining graph $\Gamma$ and use this to determine when two points in
an asymptotic cone of $A_{\Gamma}$ are separated by a cut-point.

We propose a modified version of Gersten's question: are there CAT(0)
spaces with isolated flats (\cite{hruska}) and super-linear and
sub-exponential divergence of geodesics. Our examples, like those of
Gersten and M. Kapovich do not have isolated flats. So the aphorism
may yet hold for CAT(0) spaces with isolated flats.

The organization of the paper is as follows.  In Section \ref{detour}
we recall the definitions and results concerning divergence and CAT(0)
spaces, that are pertinent to our proofs. When studying the
divergence, we use the language and techniques of detour functions
developed in \cite{macuraDET}, since they facilitate simple and
intuitive arguments.  The equivalence class of detour functions of a
proper metric space $X$ is the divergence in Gersten's sense, and it
is a quasi-isometry invariant if $X$ has a weak form of geodesic
extension property.  In Section \ref{sec:geometry} we define the
complexes $X_d$, and analyze geometric properties of the complexes
$X_d$ and $\wtd{X}_d$ pertinent to proof of the polynomial divergence
in Sections \ref{sec:proof} and \ref{sec:lower}.

In Section \ref{sec:proof} we show that the detour function of
$\wtd{X}_d$ is bounded above by a polynomial of degree $d$ and in
Section \ref{sec:lower} we show that there are geodesics $\gamma_{0}$
and $\gamma_d$ in $\wtd{X}_d$ which actually do diverge polynomially
with degree $d$, therefore establishing that the divergence of
$\wtd{X}_d$ is polynomial of degree $d$.

I am grateful to Daniel Allcock for his help with the first version of
this paper, and to the anonymous referee for careful reading and
helpful suggestions on the exposition of the paper.

\section{Detour functions and CAT(0) spaces}\label{detour}

\subsection{Detour functions and divergence} 

Detour functions were introduced in \cite{macuraDET} in order to
classify mapping tori of polynomially growing automorphisms of free
groups; they provide a language and techniques to study divergence of
geodesics in proper metric spaces, and are invariant under
quasi-isometries. We recall the definition and the main results used
in this paper, and refer the reader to \cite{macuraDET} and
\cite{gerstendiv2} for detailed expositions on detour functions and
divergence.

Let $X$ be a proper metric space, $O$ a point in $X$, and $r\geq 0 $ a
real number. Let $S(O,r)$ and $B(O,r),$ be the sphere and the open
ball, of radius $r$ centered at $O.$ We say that a path $\alpha$ in
$X$ is an $r$-detour path if $\alpha$ does not intersect $B(O, r).$
The $r$-detour distance $\delta _{r}(P,Q)$, between two points $P,Q\in
X\backslash B(O, r)$ is the infimum of the lengths of all $r$-detour
paths $\alpha$ that connect $P$ and $Q$. In the case $P$ and $Q$ are
in different components of $ X\backslash B(O, r)$, we define their
detour distance to be infinite. Since $X$ is a proper metric space, if
the detour distance $\delta _{r}(P,Q)$ is finite, Arzela-Ascoli
Theorem implies the existence of a detour path $\alpha$ such that
$|\alpha|=\delta _{r}(P,Q)$. We call such $\alpha$ a minimal length or
shortest detour path.  As indicated above, we suppress the point $O$
from notation if it is understood from the context and, when
necessary, we will talk about $(O,r)$-detour path and $(O,r)$-detour
distance.  A detour function roughly speaking, assigns to each
positive real number $r$ the maximum of all $r'\mbox{-detour}$
distances between points on the sphere of radius $r,$ where $r'=ar-b,$
for $0\leq a\leq 1$ and $b>0.$ The following definition formalizes the
above discussion:
%%%
\bdf Let $(X, d)$ be a proper metric space. Given a point $O\in X$ let
$ B_r=B(O; r) $ be the open ball of radius $r$ centered at $O$ and
$S_r=S(O;r) $ the sphere of the radius $r.$ A {\em detour function} of
$(X,O)$ is a pair $( \phi ,\mu )$ such that $\phi $ is a linear
function, $\phi (x) =ax-b,\;0<a,b;\, a\leq 1$, and $\mu :{\Bbb
  R^+}\lar {\Bbb R^+}\cup\{\infty\} $ is defined in the following way:
\[\mu_O (r) = max\{\delta _{\phi(r)}(P,Q):P,Q\in S_r\} .\]
\edf

In \cite{macuraDET} we introduced a (weak) version of geodesic
extension property for a metric space $X$ that implies the existence
of a detour function $( \phi ,\mu )$, $\phi (x) =ax-b,$ such that, if
$( \psi ,\mu' )$, $\psi (x) =cx-d,$ is a detour function and $a\leq c,
d\geq b$, then $( \phi ,\mu )$ and $( \psi ,\mu' )$, are equivalent in
the following sense. We say that $f\preceq g$ if there are constants
$A, B, C,D,E >0$ such that
$$f(x)\leq Ag(Bx+C)+Dx+E\; \mbox{for every}\; x>0.$$ We define two
functions $f,g:R^+\longrightarrow R^+\cup\{\infty\}$, to be {\em
equivalent}, $f\sim g$, if $f\preceq g$ and $g\preceq f$. This gives
equivalence relation capturing the qualitative agreement of growth
rates.  The square complexes constructed in Section \ref{sec:geometry}
satisfy a strong version of the geodesic extension property, that is
that every geodesic can be extended to an infinite geodesic ray, which
implies that any two detour functions are equivalent. In particular,
the equivalence class of the detour function does not depend on the
choice of the base point.
%%%
In the remainder of the paper, we select $a=1$ and $b=0$, and a point
$O\in X$, and take {\em the detour function} of a proper metric space
$X$ to be the function $\mu_X:{\Bbb R^+}\lar {\Bbb R^+}\cup\{\infty\}$
defined by \linebreak {$\mu_{X}(r) =\max\{\delta _{r}(P,Q):P,Q\in
S(O,r)\} .$}
%The {\em detour function} $\mu_X:{\Bbb R^+}\lar {\Bbb
%R^+}\cup\{\infty\}$ of the space $X$ is then given by \linebreak
%$\mu_X (r)= \sup\{\mu_{O}(r): O\in X\}.$

The equivalence class of the detour function of $X$ is the divergence
in Gersten's sense (see e.g.\cite{gerstendiv2}), and bounding a detour
function $\mu_X$ from above and below by polynomials of degree $d$
shows that the divergence of geodesics is polynomial of degree $d.$

\subsection{CAT(0) spaces} 

We recall the definition of a CAT(0) space and several properties that
such a space enjoys, and refer the reader to \cite{Bridson} for a
detailed treatment of the topic.

Let $(X,d)$ be a metric space and let $(E,d_E)$ be the Euclidean
plane,where $d$ and $d_E$ are the respective metrics.  A geodesic
triangle $\Delta=\Delta(P,Q,R)$ in $X$ consists of three points $P,Q,R
\in X$, its vertices, and a choice of three geodesic segments
$\gamma_{PQ}, \gamma_{QR}$ and $\gamma_{PR}$ joining the vertices, its
sides. If the point $T$ lies in the union of $\gamma_{PQ},
\gamma_{QR}$ and $\gamma_{PR}$, then we write $T\in \Delta$.

A geodesic triangle $\Delta_E=\Delta(P_E,Q_E,R_E)$ in $E$ is called a
comparison triangle for the triangle $\Delta(P,Q,R)$ if
$d(P,Q)=d_E(P_E,Q_E)$, $d(P,R)=d_E(P_E,R_E)$ and
$d(Q,R)=d_E(Q_E,R_E)$. A point $T_E$ on $\gamma_{P_EQ_E}$ is called a
comparison point for $T$ in $\gamma_{PQ}$ if $d(P,T)=d_E(P_E,T_E)$.
Comparison points for points on $\gamma_{PR}$ and $\gamma_{QR}$ are
defined in the same way.

\bdf A metric space $X$ is a CAT(0) space if it is a geodesic metric
space all of whose triangles satisfy the CAT(0) inequality:

Let $\Delta$ be a geodesic triangle in $X$ and let $\Delta_E$ be a
comparison triangle in the Euclidean plane $E$.  Then, $\Delta$ is
said to satisfy the CAT(0) inequality if for all $S,T\in \Delta$ and
all comparison points $S_E,T_E\in \Delta_E$, $d(S,T)\leq
d_E(S_E,T_E)$.\edf

A metric space $X$ is said to be of non-positive curvature if it is
locally a CAT(0) space, i.e. for every $x\in X$ there exists $r_x>0$
such that the ball $B(x,r_x)$ with the induced metric, is a CAT(0)
space.

We will use the orthogonal projections onto complete, convex, subsets
of CAT(0) spaces, called {\em projections} in \cite[II.2]{Bridson}.
We review selected parts of a proposition \cite[II.2, Proposition
2.4]{Bridson}, which gives the construction of such a projection
$\pi_C:X\lar C$.
\bpr\label{projections} Let $X$ be a CAT(0) space, and let $C$ be a
convex subset which is complete in the induced metric. Then, \ben
\item for every $x\in X$, there exists a unique point $\pi(x)\in C$
  such that $d(x,\pi(X))=d(x,C)=\mbox{inf} _{y\in C}d(x,y)$;
\item if $x'$ belongs to the geodesic segment connecting $x$ and
  $\pi(x)$, then $\pi(x)=\pi(x')$. 
\item the map $x \mapsto \pi(x)$ is a retraction from $X$ onto $C$
  which does not increase distances.

 \een 

\epr

Throughout the rest of the paper $\pi_C:X\lar C$ will denote the
projection onto complete, convex, subset $C$ of $X$, as described in
Proposition \ref{projections}. We will also make use of the property
of a CAT(0) space that a local geodesic is a global geodesic.  We will
use the following property of projections.

\brm \label{remark} Let $O$ be a point in a CAT(0) space $\wtd{X}$,
and let $O'$ be the closest point projection of $O$ to geodesic
$\omega$. If $P$ and $Q$ are points on $\omega$ such that $Q$ is
contained in the segment of $\omega$ connecting $O'$ and $P$, then
$d(O,Q)\leq d(O,P)$.  \erm

We note that the above remark is a consequence of CAT(0) inequality
applied to the triangle $OO'P$.

As a matter of general terminology and notation, a ``path $\alpha$''
refers to both the path as a continuous function $\alpha:[0,t]\lar
\mathbb{R} $ and the image of $\alpha$ in the metric space $X$, and
$|\alpha|$ stands for the length of a path $\alpha$. We will use
$\alpha\star\beta$ to denote the path which is the concatenation of
paths $\alpha$ and $\beta$, or $\gamma_1\gamma_2$ for a geodesic which
is a concatenation of geodesics paths $\gamma_1$ and $\gamma_2$.

\section{Square complexes $\wtd{X}_d$ and their geometric properties}\label{sec:geometry}

Let $G_1$ be the group $\mathbb{Z}\oplus\mathbb{Z},$ generated by
$a_0$ and $a_1$, and let $X_1$ be the flat torus obtained by isometric
identification of the edges of the Euclidean unit square $C_1=I\times
I$ . Orient the horizontal edges of the unit square from left to
right, and the vertical ones with upward positive direction. We call
this orientation {\em torus orientation}. Denote two opposite directed
edges of the square $C_1$ by $a_0$ and the other two by $a_1.$ We will
use the same notation ($a_0,$ $a_1$) for the corresponding (directed)
loops in $X_1.$ The group $G_1$ acts properly and cocompactly by
isometries on the Euclidean plane $\wtd{X}_1$ with the quotient space
$X_1$.

We define the CAT(0) groups $G_d, d\geq 2,$ inductively, taking $G_d$
to be the HNN extension of $G_{d-1}$ that amalgamates the infinite
cyclic subgroups of $G_{d-1}$ generated by $a_{d-1}$ and $a_0$
respectively.  If we denote the stable letter in this extension by
$a_d$ the resulting group $G_d$ has a presentation
\[\{a_0, \ldots, a_d\;|\; a_0a_1=a_1a_0, \;a_i^{-1}a_{0}a_i=a_{i-1},\;\mbox{for}\; 2 \leq i\leq d\}.\]
We construct a presentation complex $X_d$ of $G_d$ by a standard
topological construction of gluing with a tube, see
\cite[II.11]{Bridson} on $X_{d-1}$. Let $C_d=I\times I$ be the
Euclidean unit square with the torus orientation. Label two opposite
directed edges by $a_d$ and identify them to obtain a cylinder (tube)
$U_d.$ The remaining two edges of $C_d$ map to loops in $U_d$, and we
label them $a_{d-1}$ and $a_{0}$ respectively.  The complex $X_d$
obtained by gluing the cylinder $U_d$ to $X_{d-1}$, with the
identification map the orientation preserving isometry prescribed by
the labeling of the edges, is a graph of spaces with one vertex and
one edge. The vertex space is $X_{d-1}$ and the edge space $S^1.$ We
will call the universal cover $\wtd{X}_d$ of $X_d$ the {\it $d$-th
  square complex.} We will refer to $d$ as the {\em height} of
$\wtd{X}_d$.
   
The resulting complex $X_d$ is a non-positively curved cube complex
(see \cite[II.11]{Bridson}) and therefore $\wtd{X}_d$ is a CAT(0) cube
complex. We note that it is not difficult to see that the ``link
condition'' (\cite[II.5]{Bridson}), is satisfied: for each vertex
$P\in X_d$ every injective loop in Lk$(P,X_d)$ has length at least
$2\pi$. The preimage of $X_{d-1}$ in $\wtd{X}_d$ consists of
infinitely many disjoint convex (hence isometrically embedded) copies
of $\wtd{X}_{d-1}.$ We call such a copy of $\wtd{X}_{d-1}$ a {\em
  vertex complex} in $\wtd{X}_d$.

Let $s_d,$ $d \geq 1$ be the line segment in $C_d$ that connects the
midpoints of the opposite edges labeled $a_d.$
%%
%The is
%bi-collared and separates $X_d$.
We also denote by $s_d$ the image of $s_d$ in $U_d$, as well as its
image in $X_d$ after the gluing. Every component $H$ of the preimage
of $s_d\subset X_d$ under the covering map is isometric to the real
line and separates $\wtd{X}_d.$ We call $H$ a {\it hyperplane} in
$\wtd{X}_d$ (even though it is just a line). Every hyperplane $H$ is
geodesic contained in a single component $ostar(H)$ of the preimage of
$Int(U_d)$ under the covering map.  Since $U_d$ is a cylinder, the
closure $star(H)$ of $ostar(H)$ is isometric to a flat strip. We call
$ostar(H)$ the {\em open star} of $H$. For each edge $e$ labeled
$\td{a}_d$ in $\wtd{X}_d$, there is a unique hyperplane $H$ that
intersects $e$, and we say that $H$ {\em corresponds} to the edge $e$.

Let $\phi_0,\phi _{d-1}:S^1\lar U_d\longrightarrow X_{d}$ be the
inclusions of $S^1$ into $X_{d}$ that wrap $S^1$ isometrically once
around $a_0$ and $a_{d-1},$ respectively.  For each hyperplane $H$
there are lifts $\td{\phi}_0,\td{\phi}_{d-1}$ of $\phi_0,\phi_{d-1}$,
such that \(\omega_0=\td{\phi}_0({\Bbb R})\subset star(H)\) and
\(\omega_{d-1}=\td{\phi}_{d-1}({\Bbb R})\subset star(H).\) We call the
bi-infinite geodesic path $\omega_0$ consisting of copies of
$\td{a}_0$ {\em the smooth trace} of the hyperplane $H.$ {\em The
  rugged trace} $\omega_{d-1}$ is the bi-infinite geodesic path
consisting of copies of $\td{a}_{d-1}.$

Let $\mathcal{H}$ be the collection of all hyperplanes in $\wtd{X}_d,$
and let $U=\cup\{ostar(H_i):H_i\in\mathcal{H}\}$. Each connected
component $V$ of $\wtd{X}_d\;\backslash \;U$ is a copy of the
universal cover of the square complex $X_{d-1}$, and we called such
$V$ a vertex complex in $\wtd{X}_d.$ Since a vertex complex $V_{d-1}$
in $\wtd{X}_d$ is a copy of a square complex of height $d-1$, we can
talk about hyperplanes and vertex complexes in $V_{d-1}$. The vertex
complexes in this case are isometric copies of $\wtd{X}_{d-2}$.  In
the same fashion, for every vertex $P$ in $\wtd{X}_d$, there is a
sequence of sub-complexes $V_i\subset \wtd{X}_d$, $1\leq i\leq d$,
such that $P\in V_1\subset V_2\subset \ldots V_i\ldots \subset
V_{d-1}\subset \wtd{X}_d$.  Each $V_i$ is a copy of an $i$-th square
complex, and we call each such $V_i$ {\em an $i$-vertex complex}, or a
vertex complex of height $i$.  A hyperplane $H_i$ in $V_i$ separates
$V_i$, but does not separate $V_{i+1}$. An edge labeled $\td{a}_i$ has
{\em the height}, denoted by $height(\td{a}_i)$, equal to $i$.

If $star(H)\cap V\neq \emptyset$ for a hyperplane $H\subset \wtd{X}_d$
and a $(d-1)$-vertex complex $V$ then, $star(H)\cap V$ is a geodesic
path, the smooth or rugged trace of $H.$ We will call this geodesic
edge path the trace of $H$ in $V$ and will say that the hyperplane $H$
and the vertex complex $V$ are {\em adjacent}.

A $d$-segment (or a $d$-line) is a geodesic segment (line) that is a
concatenation of edges labeled $\td{a}_d$. When the orientation is of
importance, we will call a finite or infinite oriented segment a {\em
  ray}. We use suffixes to indicate the endpoints of a segment or the
initial endpoint of a ray: $\omega_{PQ}$ is a segment with initial
endpoint $P$ and the terminal endpoint $Q$, $\gamma_P$ stands for a
geodesic ray issuing at $P$. We call the orientation of edges in
$\wtd{X}_d$ induced by the orientation on $X_d$ the standard edge
orientation. If $\gamma_P:[0,r]\lar \wtd{X}_d$ is a geodesic ray in
1-skeleton of $\wtd{X}_d$, then $\gamma_P$ induces a
$\gamma_P$-orientation on each edge it traces, by choosing the
positive direction to be the one of increasing values of the parameter
$t\in [0,r]$. We say that a geodesic $\gamma_P$ traces an edge $e$ in
a positive direction if the $\gamma_P$-direction on $e$ coincides with
the standard direction.

\bdf\label{positive} A geodesic segment $u$ is a {\em positive}
geodesic segment if it is contained in $1$-skeleton of $\wtd{X}_{m}$,
and if $u(t)$ traces all the edges in the positive direction.  \edf

We conclude our study of basic geometric properties of square
complexes with the geodesic extension property. Since it is easily
observed that every 1-cell in $\wtd{X}_d$ is contained in a boundary
of at least two 2-cells, $\wtd{X}_d$ has no free faces, and
Proposition 5.10, Chapter II, in \cite{Bridson}, implies that
$\wtd{X}_d$ has the geodesic extension property. We formalize the
above result in the following proposition.

\bpr \label{gep} Every non-constant geodesic $\gamma$ in $\wtd{X}_d$
can be extended to an infinite geodesic ray. \epr

\section{Polynomial divergence of geodesics}\label{sec:proof}

The following theorem is our main result.

\btm \label{main} The $d$-th square complex $\wtd{X}_d$ has degree $d$
polynomial divergence of geodesics. \etm
%%%
In this section we show that there is a degree $d$ polynomial $q_d$
such that the detour function of $\wtd{X}_{d}$ is bounded above by
$q_d,$ and we start by stating and proving a lemma.
%%%
\blm\label{lemmaupper} Let $H$ be a hyperplane in $\wtd{X}_d$, and let
$\omega$ be a trace of $H$.  Let $O,Q$ be points in $ \wtd{X}_d$ such
that $Q\in star(H)$ and $d(O,Q)=r$. If $\xi\subset star(H)$ is any
bi-infinite geodesic parallel to $\omega$ such that $d(O,\xi) \leq r$,
then there is a point $P\in \xi\cap S(O,r)$, and an $(O,r)$-detour
path $\beta$ connecting $Q$ and $P$ such that $|\beta|\leq 2r+1$.
\elm

\bP 

Let $\zeta$ be the bi-infinite geodesic parallel to $\xi$ through $Q$
and let $O'$ be the projection of the point $O$ to the infinite flat
strip $U$ bounded by $\xi$ and $\zeta$. Since the lemma trivially
holds when $Q\in \xi$, we can assume that $Q\notin \xi$. This implies
that $O'\neq Q$, since otherwise we would have $d(\xi, O)>r$. Let
$\beta_{\perp}'$ be the geodesic segment in $U$ perpendicular to $\xi$
through the point $O'$, and let $U_1$ be the component of $U\backslash
\beta_{\perp}'$ that contains $Q$. We claim that the endpoints of
$\beta_{\perp}'$ are the projections of $O$ to $\zeta$ and $\xi$.  If
$O\in U$, then $O'=O$ and the claim follows directly from the
definition of $\beta_{\perp}'$.  If $O\notin U$, then either $O'\in
\zeta$, or $O'\in \xi$. If, say, $O'=O_{\zeta}\in \zeta$, then
$O_{\xi}=\beta_{\perp}'\cap \xi$ is the closest point to $O$ on $\xi$:
$O_{\zeta}$ is the point closest to $O$ on $\zeta$, and the distance
between $\zeta$ and $\xi$ is equal to the length of $\beta_{\perp}'$.
A similar argument holds if $O'=O_{\xi}\in \xi$.
\begin{figure}[h!]\label{upper1} \centering
  \includegraphics[height=2 in,width=3.2in]{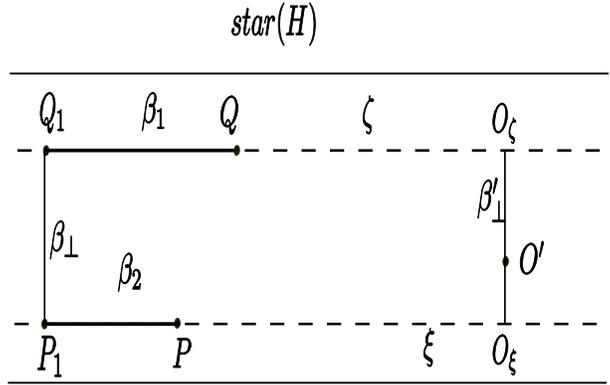}
 \caption{Illustration of the proof of Lemma \ref{lemmaupper}}
\end{figure}

Let $\beta_{\perp}$ be the geodesic segment in $U_1$ parallel to
$\beta_{\perp}'$ and at the distance $r$ from $\beta_{\perp}'$, and
denote by $Q_1$ and $P_1$ the intersections of $\zeta$ and $\xi$ with
$\beta_{\perp}$ respectively, see Figure \ref{upper1}. Since
$d(O,P_1)\geq d(O',P_1)\geq r$, and $d(O,O_{\xi})=d(O,\xi)\leq r$,
there is a point $P$ in $\xi\cap U_1$ such that $d(O,P)=r$.

To construct the desired detour path $\beta$, take $\beta_1$ to be the
segment of $\zeta$ connecting $Q$ and $Q_1$.  The length of $\beta_1$
is equal to the distance between points $Q$ and $Q_1$, which, by
construction, is not larger than $r$. Let $\beta_2$ to be the segment
of $\xi$ connecting $P_1$ and $P$, and note that the length of
$\beta_2$ is also less than or equal to $r$. Remark \ref{remark}
implies that $\beta_1$ and $\beta_2$ do not intersect $B(O,r)$, and
$\beta=\beta_1\beta_{\perp}\beta_2$ is then a detour path connecting
$Q$ and $P$ of length not larger than $2r+1$.  \eP

\blm \label{lemmabase} Let $O \in \wtd{X}_d$ and let $P,Q \in
S(O,r)\cap E$ be points contained in subcomplex $E$ of $\wtd{X}_d$
which is a copy of $\wtd{X}_1$.  Then there is an $(O,r)$-detour path
$\alpha$ connecting $P$ and $Q$ of length at most $\pi r +2r $.  \elm
%%%
\bP Let $O'$ be the projection of $O$ to $E$. By the properties of
projections, $d(O',P)\leq r$ and $d(O'Q) \leq r$. Let $\gamma_P$ and
$\gamma_Q$ be geodesics connecting $O'$ and $P,Q$ respectively.  By
the geodesic extension property we can extend $\gamma_P$ and
$\gamma_Q$ to infinite geodesic rays $\gamma_P'$ and $\gamma_Q'$.  Let
$P'$ and $Q'$ be the points of intersection of the sphere $S(O',r)$
and the geodesic rays $\gamma_P'$, $\gamma_Q'$ respectively. Let
$\beta_P$ be the segment of $\gamma_P'$, that connect $P$ and $P'$,
and let $\beta_Q$ be the segments of $\gamma_Q'$ that connect $Q$ and
$Q'$. Remark \ref{remark} implies that $d(O,S)\geq r$ for any point
$S$ in $\beta_P$. 
%if there existed a point $S'$ such that $d(O',S')\geq d(O',P) $ and
%$d(O, S') <r$, then CAT(0) comparison property would imply that
%$d(O,P)<r$, giving us a contradiction.
Similarly, $d(O,S)\geq r$ for a point $S\in \beta_Q$. Since $E$ is
copy of $\wtd{X}_1$, it is the Euclidean plane,
%geodesics in $\wtd{X}_1$ diverge
and $P'$ and $Q'$ lie on the sphere $S(O',r)$, there is an
$(O',r)$-detour path $\beta$ connecting $P'$ and $Q'$ of length at
most $\pi r$.  We note that, since the properties of projections imply
that $d(O',T)\leq d(O,T)$ for any point $T\in \beta$, $\beta$ is also
an $(O,r)$ detour path. The desired $(O,r)$-detour path that connects
$P$ and $Q$ is $\beta_P\beta\bar{\beta}_Q.$ \eP

\bpr \label{upper} There is a polynomial $q_d$, of degree $d$, such
that for any point $O$ in $\wtd{X}_d$, and any two points $P, Q$ on
the sphere $S(O,r)\subset \wtd{X}_d\,$, there is a path $\alpha$ in
$\wtd{X}_{d}\backslash B(O,r)$ connecting $P$ and $Q$ such that the
length of $\alpha$ is at most $q_{d}(r).$\epr

\bP We prove the statement of the proposition by induction.  Lemma
\ref{lemmabase} provides the base of induction, with $q_1(r)=(2+\pi)
r$.

Let $q_{d-1}$ be a polynomial of degree $d-1$ such that for any point
$O \in \wtd{X}_{d-1}$, and any two points $P', Q'\in S(O',r)\subset
\wtd{X}_{d-1}\,$ there is a path $\alpha'$ in $\wtd{X}_{d-1}\backslash
B(O,r)$ connecting $P'$ and $Q'$ such that the length $|\alpha'| \leq
q_{d-1}(r).$ Let $O, P, Q \in \wtd{X}_d$ be as in the statement of the
proposition. If $P$ and $Q$ are contained in the same vertex complex
$V_{d-1}$, the claim of the Proposition follows directly from the
induction hypothesis, otherwise, let $\mathcal{H}=\{H_1, H_2 \ldots ,
H_m\}$ be the collection of all hyperplanes in $\wtd{X}_d$ such that
each $H_i$ either separates $P$ and $Q$, or $\{P,Q\}\cap star(H_i)\neq
\emptyset$.

Without loss of generality we can assume that either $P\in star(H_1)$,
or $H_1$ is the hyperplane in $\mathcal{H}$ closest to $P$, and that
every hyperplane $H_i$ (for $i=2\ldots ,m-1$) separates $H_{i-1}$ and
$\{H_{i+1}, H_{i+2} \ldots , H_m\}.$ Since $d(P,Q)\leq 2r$, there are
no more than $2r$ hyperplanes separating $P$ and $Q$ and therefore
$m\leq 2r+2$.
%Note that every $H_i$ is at the distance at most $r$ from $O.$
Let $Y_{i}$ be the component of $\wtd{X}_d\backslash H_i$ that
contains $H_{i+1}$ and let $V_{i}$, $i=1, \ldots ,m-1$, be the
(unique) vertex complex contained in $Y_i$ that intersects
$star(H_i)$. Note that then $V_i$ also intersects $star(H_{i+1})$:
$H_{i+1}$ is contained in $Y_i$ and no hyperplane separates $H_i$ and
$H_{i+1}$. If $P\notin star(H_1)$, let $V_0$ be the vertex complex
containing $P$. Similarly, if $Q\notin star(H_m)$ let $V_m$ be the
vertex complex containing $Q$.  
%Let $O_i$ be the projection of $O$ to $V_i$.
\begin{figure}[h!]\label{upper2} \centering
  \includegraphics[height=2 in,width=4 in]{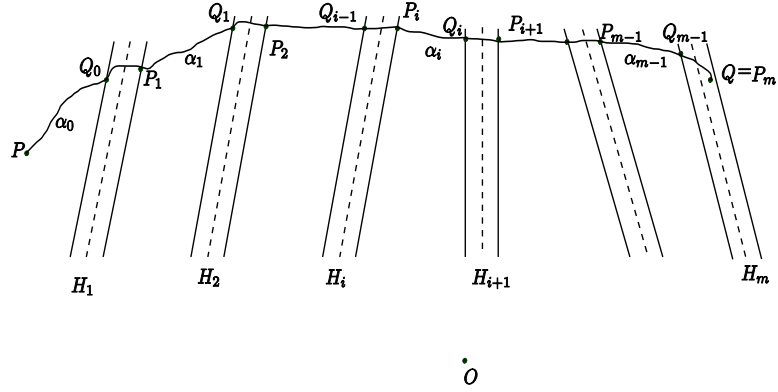}
 \caption{Illustration of the proof of Proposition \ref{upper}}
\end{figure}
If $P\in star(H_1)$ let $Q_0=P$. If $P\in V_0$, then
$d(O,star(H_1)\cap V_0)\leq r$ and we take $Q_0$ be a point in
$star(H_1)\cap V_0\cap S(0,r)$. 

If $Q\in star(H_{m})$ let $P_m=Q$, otherwise let $P_m$ be a point in
$star(H_m)\cap V_{m}\cap S(0,r)$.  By Lemma \ref{lemmaupper} there is
a point $Q_{m-1}\in V_{m-1}\cap star(H_{m})\cap S(O,r)$ and an
$(O,r)$-detour path $\beta_{m}$ connecting $Q_{m-1}$ and $P_{m}$ such
that $|\beta_{m}|\leq 2r+1$. 

For $1< i<m-2$, we let $Q_{i}$ be a point in $star(H_{i+1})\cap
V_{i}\cap S(0,r)$. Such a point exists since $d(O,star(H_{i+1})\cap
V_{i})\leq r$. By Lemma \ref{lemmaupper} for every $1\leq i<m$, there
is a point $P_{i}\in V_{i}\cap star(H_{i})\cap S(O_{i},r)$, and an
$(O,r)$-detour path $\beta_{i}$ connecting $Q_{i-1}$ and $P_{i}$ such
that $|\beta_{i}|\leq 2r+1$.
 
The point $P_i,Q_i$, chosen as above for $i=1\ldots , m-1$, are both
contained in $V_i$, and, by the induction hypothesis, for each
$i=1\ldots , m-1$ there is a detour path $\alpha_i$ of length at most
$q_{d-1}(r) $ connecting $P_i$ and $Q_i$ in the vertex space $V_i$,
and outside the ball $B(O,r)$ .

If $P\notin star(H_1)$, let $\alpha_0$ be the detour path of length at
most $q_{d-1}(r) $ connecting $P$ and $Q_0$ in the vertex space $V_0$.
Similarly, if $Q\notin star(H_m)$, let $\alpha_m$ be a detour path of
length at most $q_{d-1}(r) $ connecting $P_m$ and $Q$ in the vertex
space $V_m$.  In the case $P\in star(H_1)$ ($Q\in star(H_m)$ ) we will
take $\alpha_0$ ($\alpha_m$), to be the empty paths.

Then the path $$\alpha
=\beta_P\star\alpha_0\star\beta_1\star\alpha_1\star\beta_2\star \ldots
\alpha_{m-1}\star\beta_m\star\alpha_m\star\bar{\beta}_Q$$ is a detour
path connecting $P$ and $Q$ and $|\alpha|\leq (2r+3)q_{d-1}(r) +
(2r+2)(2r+1)$.
\eP

\section{Lower bound on the detour function}\label{sec:lower}   

We complete our proof of degree $d$ polynomial divergence in complexes
$\wtd{X}_d $ by showing that there are two geodesic rays in
$\wtd{X}_d$, emanating from the same point $O$, that diverge at least
polynomially with degree $d$.  The two such infinite rays are
$\gamma_0$ and $\gamma_d$ which are the infinite concatenations of
edges $\td{a}_0$ and $\td{a}_d$ respectively. As a matter of
convention, we use $\gamma_d$ and $\omega_d$ to denote either a
segment, a ray, or a line which is a concatenation of edges
$\td{a}_d$, and we call them $d$-segment, $d$-ray and $d$-line
respectively. We will also consider a finite oriented segment to be a
ray, issuing from the its initial endpoint.

%%%
\bdf We call the pair of geodesic rays $\gamma_0$ and $\gamma_d$ both
issuing from a vertex $T\in \wtd{X}_d$ {\em a basic $d$-corner} at $T$
and denote it by $(\gamma_d,\gamma_0)_T$.
%We call the path $\gamma_0^{-1}\gamma_d$ a basic $d$-corner path.
\edf

%%%

\bdf Let $\gamma$ and $\gamma'$ be geodesic rays in $\wtd{X}_d$. An
$(r,O)$-detour path between geodesic rays $\gamma$ and $\gamma'$ is
any $(r,O)$-detour path connecting $P\in\gamma$ and $Q\in \gamma'$,
$P,Q$ outside $B(O,r)$.  \edf
%%%
\btm\label{low} There is a polynomial $p_d$ of degree $d$ and with a
positive leading coefficient, such that the length of any
$(r,O)$-detour path in $\wtd{X}_d$ over a basic $d$-corner
$(\gamma_d,\gamma_0)_O$ is bounded below by $p_d(r)$.\etm
%%%
\subsection{Intuitive approach}

Our general approach is to prove Proposition \ref{low} by induction on
$d$. We first discuss the motivation for this approach and explain a
technical difficulty that it encounters. We start with the observation
that any detour path $\alpha\subseteq \wtd{X}_d$ over a basic
$d$-corner $(\gamma_d,\gamma_0)_O$ has to intersect every hyperplane
that $\gamma_d$ intersects. Let $n$ be the greatest integer less than
or equal to $r$, let $j$ be an integer $j\in \{1, \ldots, n\}$, and
let $H_j$ be the hyperplane that intersects $\gamma_d$ at distance
$j-1/2$ from $O$. Note that the $j$-th vertex of $\gamma_d$ (the
vertex at the distance $j$ from $O$) is contained in both $star(H_j)$
and $star(H_{j+1})$. We denote this vertex by $T_j$.  For every $1\leq
j\leq n-1$, let $\alpha_j$ be a component of $\alpha\backslash
(ostar(H_j)\cup ostar(H_{j+1}))$ that connects the rugged trace
$\omega_j$ of $H_j$ and the smooth trace $\gamma_0^{j+1}$ of $H_{j+1}$.
%$ZO_j$ and $ Z_j$ be the components of $\wtd{X}_d\backslash
%ostar(H_j)$, where $ZO_j$ is the component that contains $O$. Let
%$Y_j=Z_j\cap ZO_{j+1}$ and let $\alpha_j=\alpha\cap Y_j$.  Note that
%for every $j\leq n-1$, the path $\alpha_j$ connects the rugged trace
%$\omega_{d-1}$ of $H_j$ and the smooth trace $\omega_{0}$ of
%$H_{j+1}$ in $Y_j$.
The geodesics $\omega_{j}$, which is a $(d-1)$-ray, and
$\gamma^{j+1}_{0}$, intersect at $T_j$ and form a basic $(d-1)$-corner
at $T_j$. We note that, since $d(\alpha_j, O)\geq r$ and $d(O,T_j)=j$,
the path $\alpha_j$ does not intersect the ball of radius $r-j$
centered at $T_j$, making $\alpha_j$ into an $(r-j)$-detour over a
basic $(d-1)$-corner. We would like to use the hypothesis of induction
and claim that $|\alpha_j|\geq p_{d-1}(r-j)$, but $\alpha_j$ might not
be contained in the vertex complex $V_{d-1}$, (a copy $\wtd{X}_{d-1}$)
that contains $T_j$.  To continue the proof by induction, we would
need the hypothesis of the induction to be that the length of an
$r$-detour path over a $(d-1)$-corner in $\wtd{X}_d$ is bounded below
by $p_{d-1}(r)$. This assumption is more general than our original
statement, which brings the following additional technical difficulty
to the proof.  If $\alpha_T\subset \wtd{X}_d$ is a detour path over a
$(d-1)$-corner $(\gamma_{d-1}, \gamma_0)_T$, where $T$ is contained in
a vertex complex $V_{d-1}$, and if we do not require that
$\alpha_T\subset V_{d-1}$, then $\alpha_T$ does not necessarily
intersect the hyperplanes in $V_{d-1}$ that separate its endpoints,
making such detour paths unsuitable for induction process. We tackle
this difficulty by reformulating our statement in terms of almost
detour paths (to be defined).

The motivation for our approach is to describe a canonical way to
modify an $r$-detour path $\alpha_T$, as above, to obtain a path
$\alpha'_T$, of length not more than the length of $\alpha_T$, and
such that $\alpha'_T$ intersects all the hyperplanes in $V_{d-1}$ that
separate its endpoints. If we can then show that there is a polynomial
$p'_{d-1}$ of degree $d-1$ such that the length of $\alpha'_T$ is
bounded below by $p'_{d-1}(r)$, then $p'_{d-1}(r)$ would also give a
lower bound for the length of $\alpha_T$. The first natural question
to consider is if there is a polynomial $p'_{d-1}$ of degree $d-1$
such the length of the closest point projection of $\alpha_T$ to the
vertex complex $V_{d-1}$ is bounded below by $p'_{d-1}(r)$. We note
that $\pi_{d-1}(\alpha_T)$ might not be a detour path, and could
intersect the ball $B(O,r)$ in a collection of segments that are
rugged or smooth sides of hyperplanes in $\wtd{X}_d$.  It turns out
that the projections of detour paths are not, in general, long enough,
but if we only allow projections in the cases when they are contained
in the rugged sides of hyperplanes, we get the desired lower bound. We
proceed with this approach since, as we will show, this is sufficient
to obtain paths that behave well under induction.

In the following two subsections we introduce the terminology
necessary to define {\em almost detour paths}, which will be paths
that connect two points on the sphere $S(O,r)$, and intersect the ball
$B(O,r)$ only in geodesics of a very particular form, we will call
such geodesics {\em legal shortcuts}. Every detour path is an almost
detour path, but we will show that, given a detour path $\alpha$, we
can obtain an almost detour path $\alpha'$, with the same endpoints as
$\alpha$, and of length no longer than the length of $\alpha$, and
which has the following property: if the closest point projection
$\pi_k(\alpha'')$ of an arc $\alpha''\subset \alpha'$ to a vertex
complex $V_{k}$, $k\geq 1$ and $O\in V_k$, is a $k$-segment
$\sigma_k$, then $\pi_k(\alpha'')=\alpha''$.  This property is stated
and proved in Lemma \ref{project}, which is the most technical part of
this section. Our modified approach will also require us to consider
more general corners in addition to the basic ones, and we introduce
{\em raised corners} in the next subsection.

\subsection{Raising rays and lines} 

We recall that a positive ray is geodesic ray in $1$-skeleton of
$\wtd{X}_{m}$ that traces all its edges in a positive direction
(\ref{positive}).

\bdf A {\em raising $d$-ray $\zeta_d$} in $\wtd{X}_{m}$, $m\geq d$, is
a concatenation $\sigma_{d}u$ of a $d$-segment $\sigma_{d}$ and a
positive geodesic segment $u$, such that, if $0\leq t_1\leq t_2$, and
if $u(t_1),$ $u(t_2)$ are contained in the interiors of the edges
$e_1$ and $e_2$ respectively, then $d+1 \leq height( e_1)\leq
height(e_2)$.

\noindent 
We call $\sigma_d$ the $d$-segment of $\zeta_d$. We allow for
$\sigma_d$ to be a single point, or for $u$ to be an empty path, but
not both at the same time. In the case that $u$ is an empty path,
$\sigma_d$ cannot be a single point, and is considered to be a raising
$d$-ray.\edf
Note that a raising $d$-ray can be either an infinite ray or, a finite
segment, in which case we use the term ray to emphasize the importance
of the orientation.  It follows directly from the definition, and the
group presentation, that every raising ray is a local geodesic, and
therefore a geodesic.

\begin{figure}[h!]\label{fig4}
 \centering \includegraphics[height=3
   in,width=3.5in]{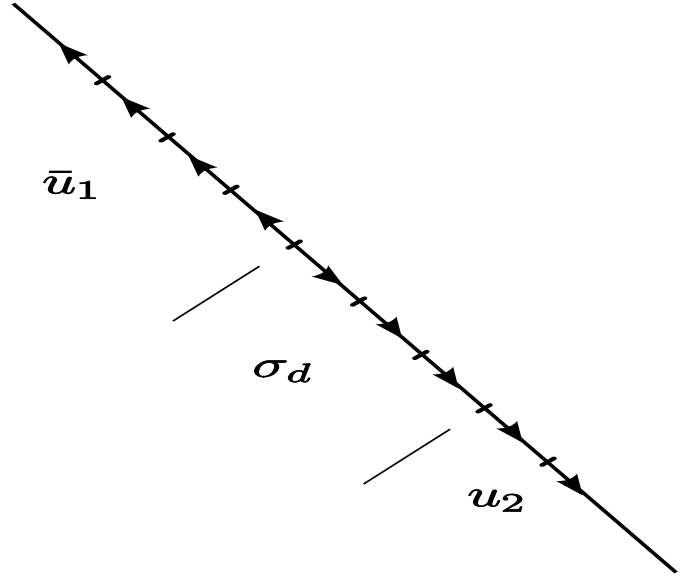}
 \caption{A raising $\sigma_d$-line}
\end{figure}

\bdf Let $\zeta_d$ be a raising $d$-ray and let $\gamma_0$ be a
$0$-segment, both issuing from a vertex $T\in \wtd{X}_m$, $m\geq
d$. We call the pair $(\zeta_d,\gamma_0)_T$ {\em a raised $d$-corner}
at $T$.  \edf
\bdf Let $\sigma_d$ be a $d$-segment in $\wtd{X}_m$, $m\geq d\geq 1$.
A geodesic $\bar{u}_1\sigma_d u_2$, where $u_1$ and $u_2$ are positive
(possibly empty) raising $(d+1)$-rays, (Figure \ref{fig4}) is called a
{\em raising $\sigma_d$-line}. \edf

\brm\label{raising} The following observation is a direct consequence
of the definition of a raising $\sigma_d$-line: if $\bar{u}_1\sigma_d
u_2$ is a raising $\sigma_d$-line and $\sigma'$ a ray contained in
$\bar{u}_1\sigma_d u_2$ such that the $\sigma'$ traces its first edge
$e$ in the positive direction, then all the edges in $\sigma'$ have
positive direction and height bigger than or equal to the height of
$e$.
%That is, $\sigma'$ is a positive raising ray.
\erm
We also note that any subray of a raising $\sigma_d$-line issuing from
a point $P$ in $\sigma_d$ is a raising $d$-ray.

\brm \label{hyper} Let $V_i\subseteq \wtd{X}_m$ be a vertex complex in
$\wtd{X}_m$, $i\leq m$, let $H_{i}$ be a hyperplane in $V_i$, and let
$\omega$ be a geodesic in $1$-skeleton of $\td{X}_m$. If $\omega$
intersects $ostar(H_i)$, then the intersection of $star(H_i)$ and
$\omega$ is contained in a single edge labeled $\td{a}_i$.\erm
\bP The only edges contained in the $ostar(H_i)$ are edges labeled
$\td{a}_i$, and, since $ostar(H_i)\cap \omega\neq \emptyset$, we
conclude that there is an edge $e$ labeled $\td{a}_i$, and a point
$Q\in Int(e)\cap \omega$. If there is a point $Q'\in star(H_i)\cap
\omega$ not contained in $e$, then convexity of $star(H_i)$ implies
that the geodesic connecting $Q$ and $Q'$ is contained in
$star(H_i)$. Moreover, since $star(H_i)$ embeds isometrically into
$\wtd{X}_m,$ such a geodesic would intersect the interior of one of
the two cubes adjacent to $e$, which contradicts the fact that the
segment of $\omega$ connecting $Q$ and $Q'$ is the unique geodesic
that connects these two points, and is contained in the $1$-skeleton
of $\wtd{X}_m$.

\eP

The above remark implies that, if $\omega$ also intersects a component
$C$ of $V_i\backslash ostar(H_i)$, then it intersects the trace of
$H_i$ adjacent to $C$ in a single point.  

\blm\label{shortcuts} Let $V_i$ be an $i$-vertex complex in
$\wtd{X}_m$, $1\leq i\leq m$, and let $H_{i}$ be a hyperplane in $V_i$
with the rugged side $\sigma_{i-1}$. If $\omega$ is a raising
$\sigma_d$-line, $1\leq d\leq i$, such that $\omega$ intersects
$\sigma_{i-1}$ at a point $P$, and such that $ostar(H_i)\cap \omega
\neq \emptyset$, then $\omega\backslash ostar(H_i)$ has exactly one
component $\omega_P$ containing $P$, and $\omega_P$ is a positive
raising $i$-ray issuing at $P$.

If $S$ is any point on $\sigma_{i-1}$, and $\sigma$ the segment of
$\sigma_{i-1}$ connecting $S$ and $P$ then the concatenation
$\sigma\omega_P$ is a raising $(i-1)$-ray.  \elm
\bP Since $\omega$ is a geodesic contained in 1-skeleton, $\omega\cap
star(H_i)$ is contained in a single edge labeled $\td{a}_i$ (Remark
\ref{hyper}). Let $Q$ be a point in $ostar(H)\cap \omega$, and let
$\omega_Q$ be the subray of $\omega$, issuing at $Q$, and that
contains $P$.  Since $P$ is contained in the rugged trace
$\sigma_{i-1}$ of $H_i$, $\omega_Q$ traces $\td{a}_i$ in the positive
direction, and Remark \ref{raising} implies that it is a positive
raising ray. Then $\omega_P=\omega_Q\backslash ostar(H_i)$ is the
component of $\omega\backslash ostar(H_i)$ containing $P$, and, since
it is a subray of $\omega_Q$, is also a positive raising ray.

The last statement of the lemma follows directly from the definition
of a raising ray.\eP

\subsection{Legal shortcuts and almost detour paths}

\bdf A {\em shortcut} is a geodesic contained in the open ball
$B(O,r)\subseteq \wtd{X}_m$.  If $\omega$ is a geodesic such that
$S(O,r)\cap \omega \neq \emptyset$ we call a point $P\in S(O,r)\cap
\omega $ an endpoint of the shortcut $\omega\cap B(O,r)$.  \edf

\bdf\label{legalshort} Let $O$ be a point in $\wtd{X}_m$. A shortcut
$\omega\subseteq \wtd{X}_m$ is {\em $O$-legal} if it is raising
$\sigma_d$-line for a $d$-segment $\sigma_d \subseteq V_d$, where
$V_d\subseteq \wtd{X}_m$ is a $d$-vertex complex, $1\leq d\leq m$,
such that $O\in V_d$.

%A $\sigma_d$-ray is {\em $d$-legal} if the $d$-segment $\sigma_d$ is
%contained in $V_d$.  
\edf 

The following properties of legal shortcuts are direct consequences of
the above definition, and we list them to provide the reader with
different aspects of legal shortcuts that we use in our proofs. \ben
\item If $\omega\subseteq \wtd{X}_m$ is an $O$-legal shortcut, and if
  $V_i\cap\omega\neq \emptyset$ for a vertex complex $V_i\subseteq
  \wtd{X}_m$, containing the point $O$, then $V_i\cap \omega$ is also
  an $O$-legal shortcut.
\item If $V_i\cap\omega= \emptyset$ for a vertex complex $V_i\subseteq
  \wtd{X}_m$, containing the point $O$, then all the edges in $\omega$
  have the height greater or equal to $i+1$.
\item If a raising $\sigma_d$-line $\omega$ for a $d$-segment
  $\sigma_d \subseteq V_d$ is an $O$-legal shortcut, and $O\in
  V_d\subset V_i$ are vertex complexes in $\wtd{X}_m$ containing $O$,
  then $\omega \cap V_i\backslash V_d$, is either empty, or consist of
  one or two positive raising $i$-rays.  \een
  
\bdf A path $\alpha\subseteq \wtd{X}_m$ is {\em an almost
  $(r,O)$-detour path} if the intersection $\alpha\cap B(O,r)$ is a
collection of $O$-legal shortcuts.  \edf
\bdf An {\em almost $(r,O)$-detour path $\alpha\subseteq \wtd{X}_m $
  over a raised corner $(\zeta_d,\gamma_0)_T$} is an almost
$(r,O)$-detour path with initial endpoint $P\in \zeta_d$, $P\notin
B(O,r)$, and such that $\alpha$ intersects $\gamma_0$ at a point
$Q\neq O$.

We also require that, if there is a shortcut $\omega\subseteq \alpha$
that contains both points $O$ and $P$, then $\omega$ is a raising
$\sigma_d$-line for a $d$-segment $\sigma_d$ (that is, has no edges of
height less than $d$).

We call the arc $\alpha'$ of $\alpha$ connecting the points $P$ and
$Q$ a {\em truncated} almost detour path. \edf

Note that we allow for $\alpha$ to contain $O$, and therefore
intersect $\gamma_0$ multiple times.

\subsection{The main result} 

\bpr\label{almostlow} For every $d\in \mathbb{N}$ there is a
polynomial $p_d$ of degree $d$, and with a positive leading
coefficient, such that for any $m,i$ such that $m\geq i\geq d$, and
any almost $(r,O)$-detour path $\alpha$ in $\wtd{X}_m$ over a raised
$i$-corner $(\zeta_i,\gamma_0)_O$, the length of the corresponding
truncated almost detour path $\alpha'$ is bounded below by
$p_d(r)$.\epr
\brm Our definition of an almost detour path over a raised corner
implies that $p_d(r)\leq p_d(r')$ for $r\leq r'$.  \erm

Since a minimal length $r$-detour path $\alpha$ over a basic
$d$-corner $(\gamma_d,\gamma_0)_O$ is also an almost $r$-detour path
over $(\gamma_d,\gamma_0)_O$, the statement of Theorem \ref{low}
follows directly from Proposition \ref{almostlow}. It remains to prove
Proposition \ref{almostlow}.

\subsection{Proof of Proposition \ref{almostlow}}

We will prove the proposition by induction on the degree of the
polynomial $p_d$ that gives the lower bound on the divergence. We
first introduce the necessary terminology, applied throughout the
statements and proofs of four lemmas that follow, and conclude the
section with the proof of Proposition \ref{almostlow}.

Let $\alpha$ be a minimal length almost $(r,O)$-detour path in
$\wtd{X}_m$ over a raised $i$-corner $(\zeta_i,\gamma_0)_O$.  For each
$k\in \mathbb{N}$, $i\leq k\leq m$, let $V_k$ be the $k$-vertex
complex in $\wtd{X}_m$ containing $O$.  Let $\pi_k:\wtd{X}_m\lar V_k$
be the projection onto $V_k$.
Let $\mathcal{E}=\{e_1, \ldots e_n\}$ be the set of all edges in
$\wtd{X}_m$ such that $\zeta_i\cap Int\,(e)\neq \emptyset$, and note
that $n\geq r$.  Since $\zeta_i$ is a geodesic, we can assume that the
distance $d(O,e_j)=j-1$ for an edge $e_j$ in $ \mathcal{E}$.

Since $\zeta_i$ is a raising ray, the $height(e_p)\leq height(e_j)$
for edges $e_p,e_j$ in $ \mathcal{E}$ such that $p\leq j$, and any
edge $e_j$ labeled $\td{a}_k$ is contained in $V_k$.  Let $H_j$ be the
hyperplane in $V_k$ that corresponds to such an $e_j$.

Then $e_j=\pi_k(e_j)$ is contained in $\pi_k(\zeta_i)$ and the
hyperplane $H_j$ separates $O$ and $\pi_k(P)$, and therefore also $Q$
and $\pi_k(P)$.  Since $Q$ and $ \pi_k(P)$ are the endpoints of
$\pi_k(\alpha)$, $H_j$ intersects $\pi_k(\alpha)$. The following lemma
establishes that, if $\alpha$ is a minimal length almost detour path,
as assumed, then $\alpha$ intersects each $H_j$.

We remark that the lemma describes the procedure in which, by
replacing arcs of $\alpha$ that do not intersect the hyperplanes by
shortcuts, any detour path $\alpha$, over a raised $i$-corner, can be
modified without increasing length to an almost detour path over a
raised $i$-corner that intersects the appropriate hyperplanes. In
particular, all the edges in the newly introduced shortcuts have
height bigger than or equal to $i$, ensuring that no shortcut in the
resulting almost detour path that contains both points $O$ and $P$
will contain edges of height strictly less than $i$.  By the same
argument, if $\alpha$ is an an almost detour path over a raised
$i$-corner, the resulting path is also an almost detour path over a
raised $i$-corner.
\blm \label{project} Let $\alpha$ be a minimal length $(r,O)$-detour
path in $\wtd{X}_m$ over a raised $i$-corner $(\zeta_i,\gamma_0).$

\noindent \ben
\item\label{1} If $\pi_k(\alpha)\cap ostar(H)\neq \emptyset$ for a
  hyperplane $H$ in $V_k$, then $\pi_k(\alpha)\cap ostar(H)\subseteq
  \alpha$.
\item \label{2} If $H$ is a hyperplane in $V_k$, where $k\geq i+1$,
  whose rugged side is contained in the vertex space $V_{k-1}$, and
  $Y_O$ the component of $V_k\backslash ostar(H)$ such that $O \in
  V_{k-1}\subseteq Y_O$, then $ \pi_k(\alpha)\subseteq Y_O$.  \een\elm
\bP Let $V_{k-1}\subseteq V_k$ be the $(k-1)$- and $k$-vertex complexes,
respectively, that contain the point $O$, and let $\pi_k:\wtd{X}_m\lar
V_k$, $\pi_{k-1}:\wtd{X}_m\lar V_{k-1}$ be the corresponding
projections.
We note that, since $\pi_m(\alpha)=\alpha$, Statement \ref{1} is true
for $k=m$, and will first prove that, for any $k$, Statement \ref{1}
implies Statement \ref{2}.

If $H$ is a hyperplane in $V_k$ whose rugged side $\sigma_{k-1}$ is
contained in the vertex space $V_{k-1}$, then $H$ does not separate
the endpoints $\pi_k(P)$ and $Q$ of $\pi_k(\alpha)$: by the definition
of a raised corner, every hyperplane in $V_k$, for $k\geq i+1$, that
separates $\pi_k(P)$ and $Q$ intersects the vertex complex $V_{k-1}$
in its smooth side.  Assume that $\pi_k(\alpha)$ intersects
$Y_C=V_k\backslash Y_O$, and note that this implies that $ostar(H)\cap
\pi_k(\alpha) \neq \emptyset$, and Statement \ref{1} of the lemma
further implies that $ostar(H)\cap \alpha=ostar(H)\cap
\pi_k(\alpha)\neq \emptyset$.  Let $\alpha_A$ be the connected
component of $\alpha \backslash ostar(H)$ that contains $P$, and let
$\alpha_B$ be the connected component of $\alpha \backslash ostar(H)$
that contains $Q$. We denote the endpoints of $\alpha_A$ and
$\alpha_B$ contained in $\sigma_{k-1}$ by $A$ and $B$
respectively. Let $\alpha'$ be the arc of $\alpha$ connecting $A$ and
$B$.

%and Note that $\alpha'$ does not intersect $\sigma_{k-1}$ except at
%its endpoints.
 
\begin{figure}[h!]\label{shorten}
  \centering
  \includegraphics[height=3.5in,width=3in]{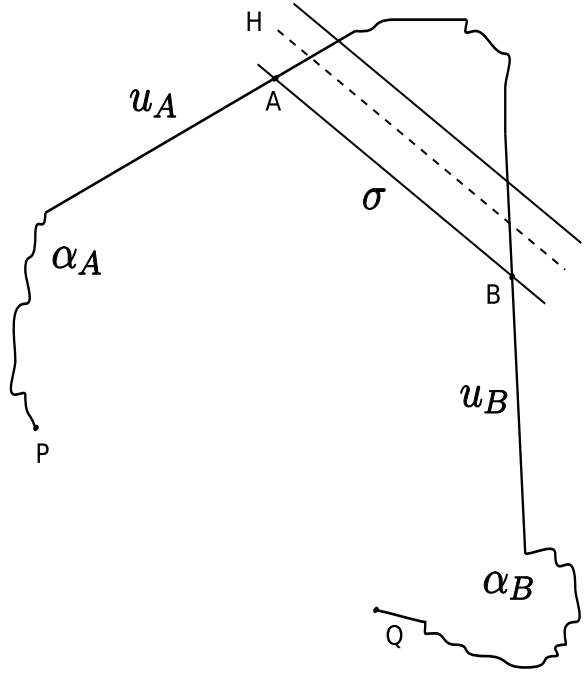}
 \caption{Shortening an almost detour path.}
\end{figure}

Let $S$ be either of the points $A, B$. If $S\in B(O,r)$, then $S$ is
contained in a legal shortcut $\omega$ which intersects $ostar(H)\cap
B(O,r)$.  Lemma \ref{shortcuts} implies that there is a positive
raising $k$-ray $u_S$ issuing from $S$, contained in $\omega$, and
such that $u_S$ does not intersect $ostar(H)$. Since $u_S$ contains
$S$ and does not intersect $ostar(H)$, $u_S\subseteq \alpha_A$. If $S$
is not in the open ball $B(O,r)$, take $u_S$ to be the trivial
path. Let $u_A$ and $u_B$ be the two positive raising rays
corresponding to the points $A$ and $B$ and let $\sigma$ be the
segment of $\sigma_{k-1}$ connecting $A$ and $B$, as in Figure
\ref{shorten}.
The path
$\theta=\bar{u}_A\sigma u_B$ is a legal shortcut with initial endpoint
in $\alpha_A$ and terminal endpoint in $\alpha_B$.  The segment
$\sigma$ is the unique geodesic segment connecting $A$ and $B$ and if
$\alpha'\neq \sigma$, or, equivalently $\alpha\cap ostar(H)\neq
\emptyset$, then $|\alpha'|>|\sigma|$. This would further imply that
$\alpha_A\sigma\alpha_B$ is an almost detour path connecting $P$ and
$Q$ of length strictly shorter than the length of $\alpha$. Since
$|\alpha_A\sigma\alpha_B|<|\alpha|$ contradicts our assumption that
$\alpha$ is an almost detour path of minimal length, we conclude that
$\alpha'=\sigma$ and therefore $\pi_k(\alpha)\subseteq Y_O$ as
claimed.

We proceed to prove that, if Statement \ref{2} of the lemma holds for
$k,$ $i+1\leq k\leq m$, then Statement \ref{1} holds for $k-1$.  First
note that the projection $\pi_{k-1}(S)$ to $V_{k-1}$ of a point
$S\notin V_{k-1}$ is contained in an edge labeled $\td{a}_{k-1}$ or
$\td{a}_0$. If such $\pi_{k-1}(S)$ is contained in $ostar(h)$ for a
hyperplane $h$ in $V_{k-1}$, then it is contained in the interior of
an edge $e$ labeled $\td{a}_{k-1}$, since edges labeled $\td{a}_{0}$
do not intersect $ostar(h)$.

Let $S$ be a point in $\alpha$ such that $\pi_{k-1}(S)$ is contained
in the interior of an edge $e\in V_{k-1}$ labeled $\td{a}_{k-1}$.
%Aiming at contradiction, assume that $S\neq\pi_{k-1}(S)$.
%%%
%We first claim that $\pi_k(S)\in V_{k-1}$ implies
%$S=\pi_k(S)=\pi_{k-1}(S)$: if $\pi_k(S)\in V_{k-1}$ then
%$\pi_{k-1}(\pi_{k}(S))=\pi_k(S)$, and $\pi_{k}(S)\in Int\,(
%e)$. Since, for a point $S\in \wtd{X}_m$ such that $\pi_k(S)\neq S$,
%$\pi_k(S)$ has to be contained in edges labeled $\td{a}_{k}$ or
%$\td{a}_0$, we conclude that $\pi_{k}(S)=S$.
%%%
Denote by $\sigma_{k-1}$ the $(k-1)$-line that contains $e$, let $H$
be the hyperplane in $V_k$ such that $ \sigma_{k-1}\subset V_{k-1}$ is
the rugged side of $H$, and let $Y_C=V_k\backslash Y_O$, where $Y_O$
is as in Statement \ref{2} of the lemma.  We assume
$S\neq\pi_{k-1}(S)$, and proceed to show that implies $\pi_k(S)\in
Y_C$, creating a contradiction to Statement \ref{2} of the lemma.  

Let $\gamma_S$ be the geodesic that connects $S$ and $\pi_{k-1}(S)$,
and let $\tau\subset Y_C$ be the smooth side of $H$.  Since $H$ is the
only hyperplane in $V_k$ such that $e\subset star(H)$, $\gamma_S$
intersects $ostar(H)$.

If $S\in ostar(H)$, then $\pi_k(S)=S \in Y_C$. If $S\notin ostar(H)$,
then $\gamma_S$ intersects $\tau$. If $S\in V_k$, then
$\gamma_S\subset V_k$, and, since $\gamma_S$ cannot intersect $\tau$
twice, $\pi_k(S)=S\in Y_C$.

If $S\notin V_k$, let $L\in \{k,\ldots,m-1\}$ be such that $S\in V_L$
and $S\notin V_{L-1}$. For $k\leq l<L$, let $H^l$ be the hyperplane in
$V_{l+1}$ adjacent to $V_{l}$, and let $\gamma^l_D$ be a geodesic in
$V_L$ connecting $S$ and a point $D$ in $V_l$. Denote by $Y_O^l$ the
component of $V_{l+1}\backslash ostar(H^l)$ such that $O\in V_l\subset
Y_O^l$, and let $Y_C^l=V_{l+1}\backslash Y_O^l$.

Since $S=\pi_L(S)\subset Y_C^{L-1}$, Statement \ref{2} implies that
$V_{L-1}$ contains the smooth trace $\tau_{L-1}$ of $H^{L-1}$, and we
also note that $\pi_{L-1}(S)\in \tau_{L-1}$. Separation properties of
hyperplanes in $V_L$ imply that any geodesic $\gamma^{L-1}_D$
intersect $\tau_{L-1}$.

\begin{figure}[h!]\label{figproject} 
 \centering
 \includegraphics[height=3.1in,width=4in]{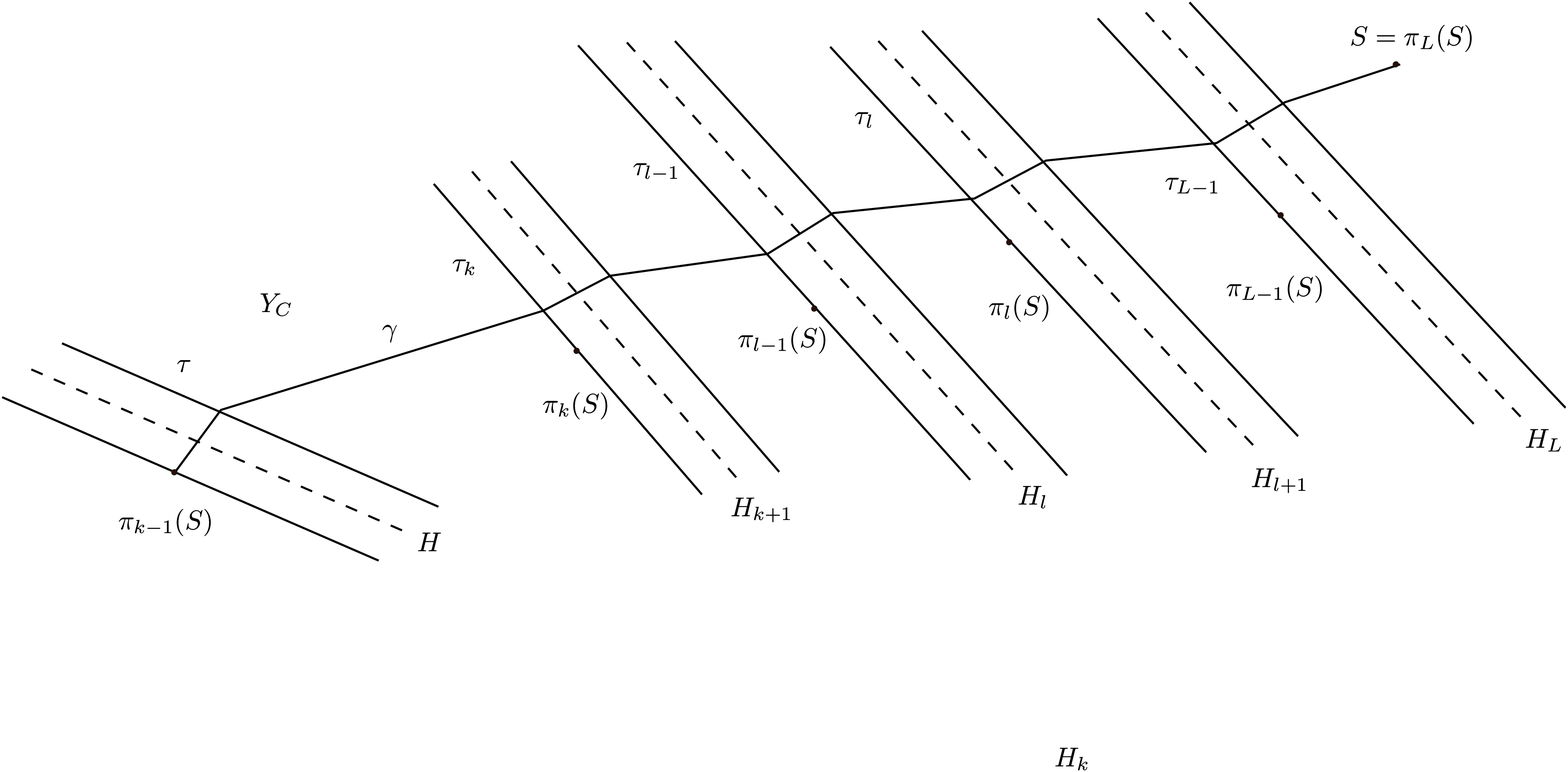}
 \caption{Projections onto vertex complexes}
\end{figure}

If we assume that $\pi_{l}(S)\in \tau_l\subset Y_C^{l-1}$, and that
any geodesic $\gamma^l_D$ intersect $\tau_{l}\subset star(H^l)$,
Statement \ref{2} implies that $V_{l-1}$ contains the smooth trace
$\tau_{l-1}$ of $H^{l-1}$, and separation properties of hyperplanes in
$V_l$ imply that any geodesic $\gamma^{l-1}_D$ intersect
$\tau_{l-1}$. We, deduce, by induction that the $V_l$ contains the
smooth trace $\tau_l$ of $H^l$ for every $k\leq l<L$, that
$\pi_l(S)\in \tau_l(S)$, and also that, for every such $l$, any
$\gamma^l_D$ intersects $\tau_l$, see Figure \ref{figproject}.
We turn our attention back to $\gamma_S$, concluding, based on the
above discussion, that it intersect $\tau_k$, and that $\tau_k$ also
contains $\pi_k(S)$. Since $\gamma_S$ cannot intersect $\tau$ twice
$\tau$ is contained in the component of $V_k\backslash ostar(H)$ that
contains $\tau_k$, which is $Y_C$, and therefore $\pi_k(S)\in Y_C$.

Since the the rugged side of $H$ is contained in $V_{k-1}$, this
provides the contradiction, which we aimed for, to Statement \ref{2}
of the lemma that $\pi_k(\alpha)\cap Y_C=\emptyset$. We conclude that
$S=\pi_{k-1}(S)$, establishing our claim that $\pi_{k-1}(\alpha)\cap
ostar(h)$ is contained in $\alpha$.  \eP

We highlight the following straightforward corollary of the second
statement of Lemma \ref{project}.
\bcr\label{cor} If $H$ is a hyperplane in $V_k$, where $k\geq i+1$,
whose rugged side is contained in the vertex space $V_{k-1}$ then
$\alpha$ does not intersect $H$. \ecr

\blm\label{noproper} We may assume that no proper arc of $\alpha$ is
an almost detour path over a raised $i$-corner, $i\geq d$, based at
$O$. \elm

\bP If there is such a proper arc $\alpha'$, we can replace $\alpha$
by $\alpha'$.  Since $\alpha$ has finite length and
$|\alpha'|<|\alpha|$, this process terminates.  \eP

In particular, Lemma \ref{noproper} implies that there are no raising
$i$-rays, $i\geq d$, issuing at $O$ and intersecting $\alpha\backslash
B(O,r)$ at any point other than $P$. We next show that such an almost
detour path does not intersect $\zeta_i$ at any other points except
$P$ and, possibly, $O$.

\blm\label{noproper2} $\alpha\cap \zeta_i\subset\{O,P\}$.  \elm

\bP If $\alpha$ intersects $\zeta_i$ in a point other then $\{O,P\}$,
then $\alpha$ contains a shortcut $\omega$ that intersects $\zeta_i$
in a point other then $\{O,P\}$.

We first consider the case that such a shortcut $\omega$ contains $O$.
Then $O$ separates $\omega$ into two components $\omega_1$ and
$\omega_2$. If one of the components, say $\omega_1$, contains $P$
then $\omega$ is an raising $i$-line, and $\omega_2$ is a raising
$i$-ray that intersects $\alpha\backslash B(O,r)$ at a point other
than $P$.  If neither of the rays $\omega_1$, $\omega_2$ contain $P$,
but if one of them, say $\omega_1$, intersects $\zeta_i$ at a point
other than $O$, then $\omega_1$ is a raising $i$-ray that intersects
$\alpha\backslash B(O,r)$ at a point other than $P$. In either case,
we have a contradiction with Lemma \ref{noproper}.
 
If $\omega$ does not contain $O$, and intersect $\zeta_i$ at point(s)
different then $P$, let $Z_0\neq P$ be the point in $\omega\cap
\zeta_i$ closest to $O$.  $Z_0$ separates $\alpha$ into arcs
$\alpha_P$ and $\alpha_Q$ containing $P$ and $Q$ respectively, and
also separates $\omega$ into two components. We denote the component
of $\omega\backslash \{Z_0\}$ contained in $\alpha_Q$ by
$\omega_0$. We consider $\omega_0$ as a ray issuing from $Z_0$.

Let $\delta$ be the segment of $\zeta_i$ connecting $O$ and $Z_0$, and
let $e_j$ be the edge in $\delta$ such that $Z_0$ is an endpoint of
$e_j$. We want to show that the height of every edge in $\omega$ is
greater than or equal to the height of $e_j$. We start by proving that
$\omega$ does not intersect $H_j$, the hyperplane that corresponds to
the edge $e_j$. By the definition of the point $Z_0$ and the segment
$\delta$, $\omega$ does not intersect $e_j$ at any point other than
$Z_0$. If $\omega$ intersects $H_j$ at a point $Z_1 \notin e_j$, then
the segment of $\omega$ connecting $Z_0$ and $Z_1$ would be the
geodesic between $Z_0$ and $Z_1$ in $star(H_j)$.  Since such a
geodesic is not contained in 1-skeleton, it cannot be a part of a
shortcut.

Let $h_j$ be the height of the edge $e_j$, and, consistent with our
notation, let $V_{h_j}$ be the $h_j$-vertex complex containing
$O$. Since $\omega$ does not intersect $H_j$, the closest point
projection $\pi_{h_j}(\omega)$ does not intersect $H_j$ either (Lemma
\ref{project}, statement (1)).

Note that $H_j$ separates $V_{h_j}$ into two components, one of which
contains $O$ and $V_{h_j-1}$, and the other one containing
$Z_0$. Since $\pi_{h_j}(\omega)$ contains $Z_0$ and does not intersect
$H_j$, it does not intersect $V_{h_j-1}$ either. This further implies
that, since $\pi_{h_j}(\omega\cap V_{h_j} )=\omega\cap V_{h_j}$,
$\omega$ does not intersect $V_{h_j-1}$.

Since $\omega$ is an $O$-legal shortcut, the observation (item 2)
after definition \ref{legalshort} implies that the height of every
edge in $\omega$ is greater or equal to the height of $e_j$, which is
greater or equal to $i$. 

If the height of the first edge $e$ that $\omega_0$ traces is equal to
the height of $e_j$, then $\delta\omega_0$ traces $e_j$ and $e$ in the
same direction. If the height of $e$ is strictly larger than the
height of $e_j$, then Corollary \ref{cor} implies that $\omega_0$
traces all its edges in positive direction. This, together with
properties of legal shortcuts, makes $\delta\omega_0$ a raising
$i$-ray that intersects $\alpha \backslash B(O,r)$ at a point
different than $P$, contradicting Lemma \ref{noproper} again. \eP

The following lemma establishes the base of induction for our
inductive proof of Proposition \ref{almostlow}. 

\blm\label{baseind} Let $\alpha'$ be a minimal length almost
$(r,O)$-detour path over a raised \mbox{$d$-corner} $(\zeta_d,
\gamma_0)$ based at $O$, with the initial endpoint $P\in \zeta_d$.
Let $Q$ be the point of intersection of $\alpha'$ and $\gamma_0$, and
let $\alpha$ be the arc of $\alpha'$ connecting $P$ and $Q$. Then the
length of $\alpha$ is at least $r-1$.  \elm

\bP The claim follows directly from the above observation that
$\alpha$ intersects both traces of $star(H)$ for every hyperplane
$H_j$, $1\leq j \leq n-1$: the distance between such two intersections
is at least 1, which implies that the length of $\alpha$ is at least
$r-1$. \eP

We can now complete the proof of Proposition \ref{almostlow}.

By Lemma \ref{baseind}, for $i\geq 1$, the length of a truncated
almost $r$-detour path over a raised $i$-corner is at least
$p_1(r)=r-1$, which establishes the base of induction.  We proceed to
prove that the existence of a polynomial $p_{d-1}$ of degree $d-1$,
where $d\geq 2$, as in the statement of the proposition, implies the
existence of a polynomial $p_{d}$.

\begin{figure}[h!]\label{fig1} \centering
  \includegraphics[height=3.5in,width=3.2in]{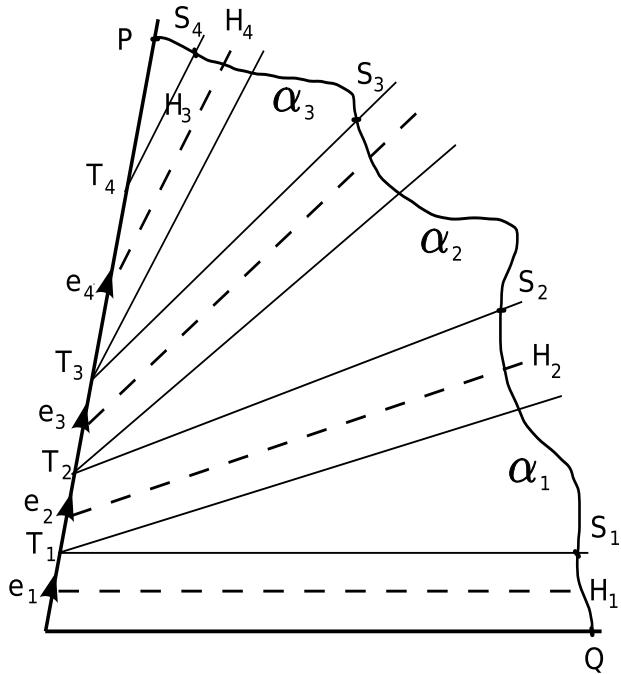}
 \caption{Intersections of hyperplanes and an almost detour path
   $\alpha$ when $\zeta_i$ traces the edges in $\sigma_i$ in positive
   direction. }
\end{figure}

Let $\zeta_i=\sigma_i u_{i+1}$, where $\sigma_i$ is an $i$-segment and
$u_{i+1}$ is a raising $(i+1)$-ray.  For each edge $e_j$ in
$\sigma_iu_{i+1}$, $j\leq n-1$, we denote by $T_j$ the terminal
endpoint of $e_j$. Since $\zeta_i$ traces each edge $e_j$ in $u_{i+1}$
in positive direction, the point $T_j$ lies in the rugged trace of the
hyperplane corresponding to $e_{j}$, and on the smooth trace of the
hyperplane corresponding to the edge $e_{j+1}$.
If $\zeta_i$ also traces the edges in $\sigma_i$ in positive
orientation, the same conclusion holds for all $T_j$ in $\{T_1, \ldots
T_{n-1}\}$: $T_j$ is contained in the rugged trace of the hyperplane
$H_j$ and the smooth side of the hyperplane $H_{j+1}$, as illustrated
in Figure \ref{fig1}.

\begin{figure}[h!]\label{fig2} 
 \centering
 \includegraphics[height=3.5in,width=3.2in]{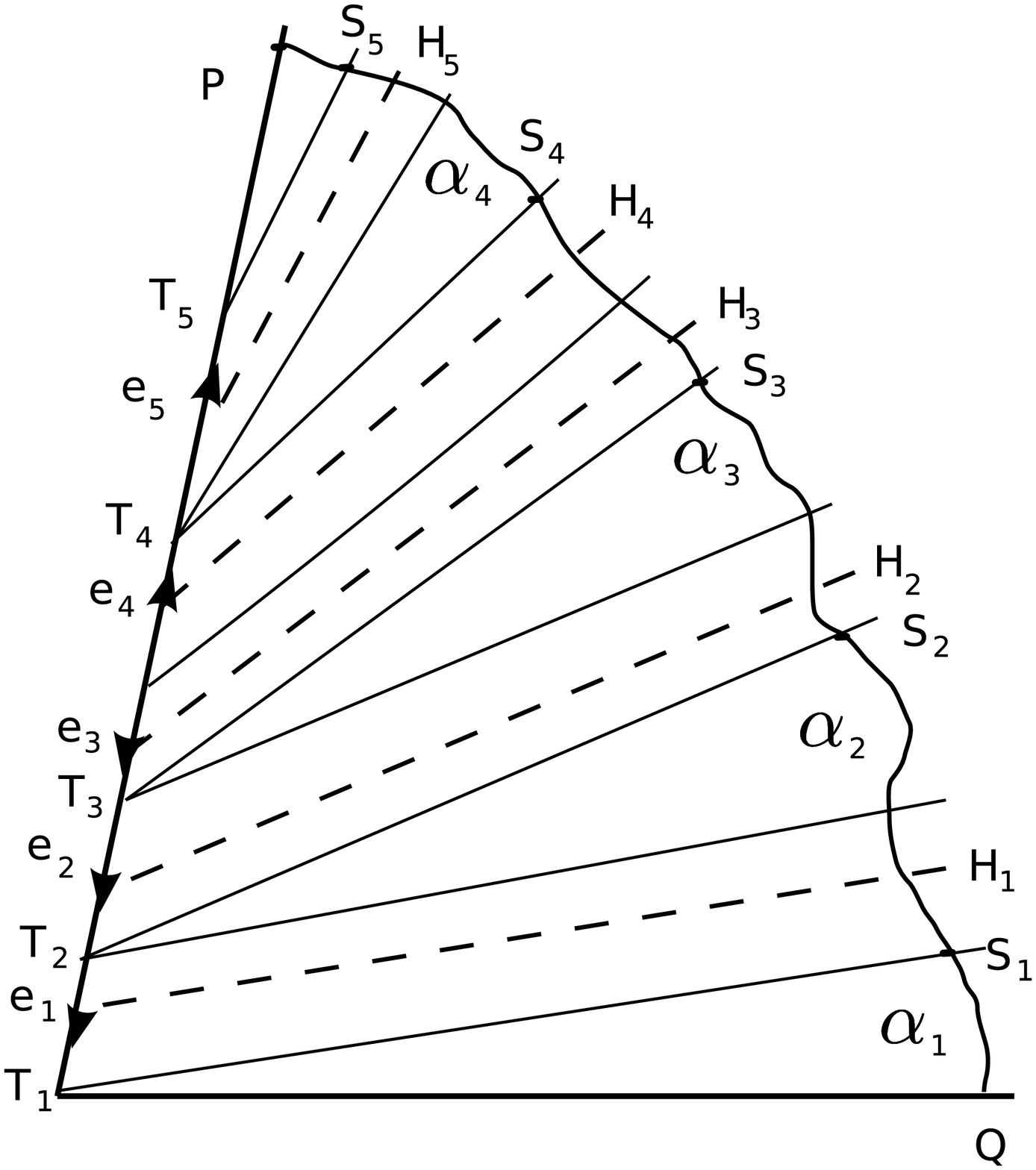}
 \caption{Intersections of hyperplane and an almost detour path
   $\alpha$ when $\zeta_i$ traces the edges in $\sigma_i$ in negative
   direction.}
\end{figure}

If $\zeta_i$ traces the edges in $\sigma_i$ in negative direction, as
illustrated in Figure \ref{fig2}, then, for $j\geq 2$, the vertex $T_j$
which is the terminal endpoint of an edge $e_j\in \sigma_i$ lies in
the rugged trace of the hyperplane corresponding to the edge $e_{j}$,
and on the smooth trace of the hyperplane corresponding to the edge
$e_{j-1}$.  The vertex $T_{1}$ lies on the rugged trace of the
hyperplane corresponding to the edge $e_{1}$ and on the smooth path
$\gamma_0$.

We will show that there is a collection $\{\alpha_j |j\in\{1,\ldots ,
n-1\}\}$ of disjoint arcs $\alpha_j$ of $\alpha$, such that $\alpha_j$
is a truncated almost detour path over a raised $(i_j-1)$-corner,
$i_j\geq d$, based at $T_j\in V_k$.

Let $q$ to be the element of $\{0,1,\ldots n-1 \}$ defined in the
following way.  If $\zeta_i$ traces all the edges $e_j$, $j\in
\{1,\ldots n-1\}$, in the positive direction, we let $q=0$.  If
$\zeta_i=\sigma_i$ and if it traces all the edges $e_j$, $j\in
\{1,\ldots n-1\}$, in the negative direction, we let $q=n-1$.
Otherwise, let $q$ be such that $\zeta_i$ traces each edge $e_j$ in
negative direction for $j\leq q$, and in positive direction the for
$j>q$.
Let $S_n=P$, and $S_0=Q$. For each $j\in \{1,\ldots n-1 \}$ we
inductively define a point $S_j$ in the rugged side of the hyperplane
$H_j$ in the following way. If $j$ is such that $0\leq q <j \leq n-1$,
let $S_j$ to be the point of intersection of $\alpha$ and the rugged
side of the hyperplane $H_j$ such that the arc $\beta_j$ of $\alpha$
connecting $S_j$ and $P$ does not intersect $H_j$.  If $1 \leq j \leq
q$ choose $S_j$ be the point in the intersection of $\alpha$ and the
rugged side of the hyperplane $H_j$ such that the arc $\beta_j$ of
$\alpha$ connecting $S_j$ and $Q$ does not intersect $H_j$. Lemma
\ref{noproper2} implies that $S_j\neq T_j$

If $S_ j$ is outside the ball $B(O,r)$ we let $P_j=S_j$.  If $S_j\in
B(O,r)$, then it is contained in a shortcut $\omega$ such that
$\omega\cap ostar(H_j)\neq \emptyset $, and, since $S_j$ is contained
in the rugged side of $H_j$, Lemma \ref{project} implies that there is
an endpoint $P_j$ of $\omega$ such that the oriented segment of
$\omega$ connecting $S_j$ and $P_j$ is an positive raising ray which
does not intersect $H_j$, and such that the geodesic connecting $T_j$
and $P_j$ is a raising $(i_j-1)$-ray. We note that $P_j$ is contained
in $\beta_j$.

Let $\zeta^j$ be the raising $(i_j-1)$-ray connecting $T_j$ and $P_j$,
and let $\gamma^j_0$ the $0$-ray issuing at $T_{j}$. We note that, for
$j>q$, $\zeta^j$ does not contain $S_{j+1}$: if it did, the segment of
$\zeta^j$ between $T_j$ and $S_{j+1}$ would be the unique geodesic
segment connecting $T_j$ and $S_{j+1}$, and, the definition of
$S_{j+1}$ together with our observation that $S_{j+1}\neq T_{j+1}$,
imply that such geodesic segment is not contained in 1-skeleton. By
the same argument, $S_{j-1}$ is not contained in $\zeta^j$ for $j\leq
q$.

We let $\alpha_j$ be the arc of $\beta_j$ connecting $P_j$ and
$P_{j+1}$ in the case $j\geq q$, and the arc of $\beta_j$ connecting
$P_j$ and $P_{j-1}$ for $j<q$. The above discussion implies that the
point $P_j$ is contained in the arc of $\beta_j$ connecting $S_j$ and
$S_{j+1}$ for $j>q$, and in the arc of $\beta_j$ connecting $S_j$ and
$S_{j-1}$ for $j\leq q$, and therefore the arcs $\alpha_j$ are
disjoint.  Each $\alpha_j$ is an almost detour path and Lemma
\ref{noproper2} implies that $\alpha_j$ intersects $\gamma_0^j$ at a
point different than $T_j$. Therefore, for every $j\in \{1,\ldots,
n-1\}$, $\alpha_{j}$ is a subarc of $\alpha$, which is a truncated
almost detour path over a $(i_j-1)$-corner $(\zeta^j, \gamma_0^j)$ and
Since $d(O,T_j)\leq j$, and $\alpha$ is an an almost $(r,O)$-detour
path, $d(T_{j},A) \geq r-j$ for every point $A\in \alpha_j$ which is
not contained in a legal shortcut.  Therefore $\alpha_j$ is an
$(r-j)$-almost detour path over a raised $i_j-1$-corner, for $i_j\geq
d$. By the hypothesis of the induction $|\alpha_j|\geq
p_{d-1}(r-j)$. Then the length of $$|\alpha|\geq \sum
_{j=1}^{n-1}|\alpha_j|\geq \sum _{j=1}^{n-1}p_{d-1}(r-j), $$ and
$p_d(r)=\sum _{j=1}^{n-1}p_{d-1}(r-j) $ is a polynomial of degree $d$
which is a lower bound for the detour function of $\wtd{X}_m$.

\bibliographystyle{plain}

\bibliography{reference}
\end{document}